\title{Output feedback stabilization of the Korteweg-de Vries equation\footnote{This work has been partially supported by Fondecyt 1140741,
and Conicyt-Basal Project FB0008 AC3E.}}

\author{Swann Marx\footnote{GIPSA-lab, Department of Automatic Control, Grenoble Campus, 11 rue des math\'ematiques, BP 46, 38402 Saint Martin d'H\`eres Cedex, France. E-mail: swann.marx@gipsa-lab.com} \: and\: Eduardo Cerpa\footnote{Departamento de Matem\'atica, Universidad T\'ecnica Federico Santa Mar\'ia, Avda. Espa\~na 1680, Valpara\'iso, Chile. E-mail: eduardo.cerpa@usm.cl}}
\documentclass[12pt]{article}

\usepackage{amsfonts,latexsym,amsmath,amssymb}
\usepackage[mathscr]{eucal}
\usepackage{graphicx,color}
\usepackage{comment}
\usepackage{hyperref}
\usepackage{epstopdf}
\newtheorem{nntheorem}{Theorem}

\newtheorem{remark}{Remark}

\catcode`\@=11
\def\downparenfill{$\m@th\braceld\leaders\vrule\hfill\bracerd$}
\def\overparen#1{\mathop{\vbox{\ialign{##\crcr\crcr
\noalign{\kern0.4ex}
\downparenfill\crcr\noalign{\kern0.4ex\nointerlineskip}
$\hfil\displaystyle{#1}\hfil$\crcr}}}\limits}
\catcode`\@=12

\begin{document}
\maketitle

\begin{abstract}                          
This paper presents an output feedback control law for the Korteweg-de Vries equation. The control design is based on the backstepping method and the introduction of an appropriate observer. The local exponential stability of the closed-loop system is proven. Some numerical simulations are shown to illustrate this theoretical result.  
\end{abstract}





\section{Introduction}\label{intro}

The Korteweg-de Vries (KdV) equation was introduced in 1895 to
describe approximatively the behavior of long waves in a water channel of relatively shallow depth. Since then, this equation has attracted a lot of attention due to fascinating mathematical features and a number of possible applications. 

From a control viewpoint, the KdV system also presents amazing behaviors. Surprisingly, 
by considering di\-ffe\-rent boundary actuators on a bounded interval $[0,L]$, we get control results of different nature. Roughly speaking, the system is exactly controllable when the control acts from the right endpoint $x=L$ (\cite{rosier1997kdv,cerpa-rivas-zhang}), and null-controllable when the control acts from the left endpoint $x=0$ (\cite{glass-guerrero,guilleron,carreno-guerrero}). 

Due to this kind of phenomena, the control properties of this nonlinear dispersive partial differential equation 
have been deeply studied. However, there still are many open questions. See \cite{cerpa2013control}, \cite{rosier-zhang}, \cite{mcpa-siam}, and the references therein.

In this article we focus on the boundary stabilizability problem for the KdV equation with a control acting on the left Dirichlet boundary condition. 
The studied system can be written as follows
\begin{equation}
\label{lKdV_equation_ld}
\left\{
\begin{split}
&u_t+u_{x}+u_{xxx}+uu_x=0,\\
&u(t,0)=\kappa(t),\: u_x(t,L)=0,\: u_{xx}(t,L)=0,\\
&u(0,x)=u_0(x),
\end{split}
\right.
\end{equation}
where $\kappa=\kappa(t)$ denotes the boundary control input and $u_0=u_0(x)$ is the initial condition. Concerning the stability when no control is applied ($\kappa=0$), it is known that the linear system is asymptotically stable (see \cite[Lemma 3]{tang2013stabKdV}). We aim here to design an output feedback control in order to get the exponential stability of the closed-loop system.





Some full state feedback controls have already been designed in the literature for KdV systems. When the control acts on the right endpoint, we find \cite{cerpa2009rapid} where a Gramian-based method is applied, and \cite{coron-lu} where some suitable integral transforms are used. In   \cite{tang2013stabKdV}, \cite{cerpa_coron_backstepping} and \cite{marxcerpa2014outputkdv}, the authors use the backstepping method to design feedback controllers acting on the left endpoint of the interval.

However, in most cases, we have no access to measure the full state of the system. Thus, it is more realistic to design an output feedback control, i.e.,  a feedback law depending only on some partial measurements of the state. 

For autonomous linear finite-dimensional systems, the separation principle holds. Thus, stabilizability and observability assumptions are sufficient to ensure the stability of the closed-loop system. In other words, if there exists a controller, which asymptotically stabilizes  the origin of the system and an observer which converges asymptotically to the state system, the output feedback built from this observer and this state feedback asymptotically stabilizes  the origin of the system. The case of autonomous nonlinear finite-dimensional systems depends critically on the structure of the system. We can only hope having a semi-global result (see e.g \cite{firststep}). In a PDE framework, this principle, even for linear systems, is no longer true and the stability of the closed-loop system is not guaranteed. 

The basic question to state the problem is which kind of measurements we are going to consider, being the boundary case  the most challenging one. In \cite{marxcerpa2014outputkdv} we consider the linear KdV equation with boundary conditions
\begin{equation}\label{ii}u(t,0)=\kappa(t),\: u(t,L)=0,\: u_{x}(t,L)=0.\end{equation}
In that paper, we see that this system is not obser\-va\-ble from the output $y(t)=u_x(t,0)$ for some values of $L$. However, we design an output feedback law exponentially stabilizing the system for the output given by $y(t)=u_{xx}(t,L)$. Thus, we see that the choice of the output is crucial. In \cite{hasan2016output}, the same  controller has been applied to the  nonlinear Korteweg-de Vries equation. 

In this paper, we will consider the nonlinear KdV equation \eqref{lKdV_equation_ld} with measurement
\begin{equation}
\label{measurement}
y(t)=u(t,L).
\end{equation}
Independently to \cite{marxcerpa2014outputkdv} and to the present paper, Tang and Krstic have developed the same program for similar linear KdV equations. Full state \cite{tang2013stabKdV} and output state \cite{output-tang-kristic} feedback controls are designed by using the backstepping method.


This paper is organized as follows. In Section \ref{main_results}, we state our main result. Section \ref{state_feedback_design} is devoted to recall the feedback control designed in \cite{cerpa_coron_backstepping}. The observer is built in Section \ref{observer_design}. In Section \ref{stability_analysis_output_feedback}, the stability of the linear closed-loop controller-observer system is proven. In Section \ref{n-l-kdv}, we prove the local stability of the nonlinear closed loop controller-observer system. Some numerical simulations are presented in Section \ref{simu}. Finally, Section \ref{con} states some conclusions.\\
\section{Main Result}
\label{main_results}

Based on \cite{krstic_smyshlyaev_backstepping} and \cite{smyshlyaev2005backstepping}, we built an observer for \eqref{lKdV_equation_ld}. More precisely, we define, for some appropriate function $p_1(x)$, the following copy of the plant with a term depending on the observation error 
\begin{equation}
\label{linear_KdV_observer}
\left\{
\begin{split}
&\hat{u}_t+\hat{u}_x+\hat{u}_{xxx}+\hat{u}\hat{u}_x+p_1(x)[y(t)-\hat{u}(t,L)]=0,\\
&\hat{u}(t,0)=\kappa(t),\: \hat{u}_x(t,L)=\hat{u}_{xx}(t,L)=0,\\
&\hat{u}(0,x)=\hat u_0.
\end{split}
\right.
\end{equation}

As mentioned in Section \ref{intro}, our main result is the local stabilization of the KdV equation by using the  output \eqref{measurement}, as stated in the following theorem whose proof is given in Section \ref{n-l-kdv}.\\

\begin{nntheorem}\label{main-th} For any $\lambda>0$, there exist an output feedback law $\kappa(t):=\kappa(\hat{u}(t,x))$, a function $p_1=p_1(x)$, and two constants $C>0$, $r>0$ such that for any initial conditions $u_0,\hat u_0 \in L^2(0,L)$ satisfying
\begin{equation}
\Vert u_0\Vert_{L^2(0,L)}\leq r ,\quad \Vert \hat u_0\Vert_{L^2(0,L)}\leq r,
\end{equation}
the solution of \eqref{lKdV_equation_ld}-\eqref{measurement}-\eqref{linear_KdV_observer} satisfies
\begin{multline}
\Vert u(t,\cdot)-\hat{u}(t,\cdot)\Vert_{L^2(0,L)}+\Vert \hat{u}(t,\cdot)\Vert_{L^2(0,L)}
\leq C e^{-\lambda t}\Big(\Vert u_0-\hat u_0\Vert_{L^2(0,L)}+\Vert \hat u_0 \Vert_{L^2(0,L)} \Big),\quad \forall t\geq 0.
\end{multline}

\end{nntheorem}

\begin{remark} Notice that from this theorem we get the exponential decreasing to $0$ of the $L^2$-norm of the solution $u=u(t,x)$ provided that the $L^2$-norm of the initial conditions of the plant and the observer are sufficiently small. 
\end{remark}

\section{Control Design}
\label{state_feedback_design}

The backstepping design applied here is based on the linear part of the equation. Thus, we consider the control system linearized around the origin
\begin{equation}
\label{lkdv}
\left\{
\begin{split}
&u_t+u_x+u_{xxx}=0,\\
&u(t,0)=\kappa(t),\: u_x(t,L)=u_{xx}(t,L)=0,
\end{split}
\right.
\end{equation} 
and the linear observer
\begin{equation}
\label{lkdv_ob}
\left\{
\begin{split}
&\hat{u}_t+\hat{u}_x+\hat{u}_{xxx}+p_1(x)[y(t)-\hat{u}_{xx}(t,L)]=0,\\
&\hat{u}(t,0)=\kappa(t),\: \hat{u}_x(t,L)=\hat{u}_{xx}(t,L)=0.
\end{split}
\right.
\end{equation}

The standard method of output feedback design follows a three-step strategy. We first design the full state feedback control. Next, we built the observer. Finally, we prove that plugging the observer state into the feedback law stabilizes the closed loop system. 



 In \cite{cerpa_coron_backstepping} the following Volterra transformation is introduced
\begin{equation}\label{mapk}
w(x)=\Pi(u(x)):=u(x)-\int_x^Lk(x,y)u(y)dy.
\end{equation}
The function $k$ is chosen such that $u=u(t,x)$, solution of (\ref{lkdv}) with control
\begin{equation}\label{feedback} \kappa(t)=\int_0^Lk(0,y)u(t,y)dy,
\end{equation}
is  mapped into the trajectory $w=w(t,x)$, solution of the linear system
\begin{equation}
\label{sf_voltransf}
\left\{
\begin{split}
&w_t+w_x+w_{xxx}+\lambda w=0,\\
&w(t,0)=w_x(t,L)=w_{xx}(t,L)=0,
\end{split}
\right.
\end{equation}
which is exponentially stable for $\lambda >0$, with a decay rate at least equal to $\lambda$.

The kernel function $k=k(x,y)$ is characterized by
\begin{equation}
\label{sf_kernel}
\left\{
\begin{split}
&k_{xxx}+k_{yyy}+k_x+k_y=-\lambda k,\text{ in }\mathcal{T},\\
&k(x,L)+k_{yy}(x,L)=0,\text{ in } [0,L],\\
&k(x,x)=0,\text{ in } [0,L],\\
&k_x(x,x)=\frac{\lambda}{3}(L-x),\text{ in }[0,L],
\end{split}
\right.
\end{equation}
where $\mathcal{T}:=\lbrace (x,y)/x\in [0,L],y\in [x,L]\rbrace$. The solution of (\ref{sf_kernel}) exists. This is proved in \cite[Section VI]{cerpa_coron_backstepping} by using the method of successive approximations. Unlikely the case of heat or wave equations, we do not have an explicit solution.

In \cite{cerpa_coron_backstepping} it is proved that the transformation \eqref{mapk} linking (\ref{lKdV_equation_ld}) and (\ref{sf_voltransf}) is invertible, continuous and with a continuous inverse function. Therefore, the exponential decay for $w$, solution of \eqref{sf_voltransf}, implies the exponential decay for the solution $u$ controlled by \eqref{feedback}. Thus, with this method, the following theorem is proven.\\

\begin{nntheorem} (\cite{cerpa_coron_backstepping})
For any $\lambda>0$, there exists $C>0$ such that
\begin{equation}
\Vert u(t,\cdot)\Vert_{L^2(0,L)}\leq Ce^{-\lambda t}\Vert u(0,\cdot)\Vert_{L^2(0,L)}, \quad \forall t\geq 0,
\end{equation}
for any solution of \eqref{lkdv}-\eqref{feedback}.
\end{nntheorem}


We give later more details on the observer design. Let us remark that the output feedback  law is designed as
\begin{equation}\label{kk-obs}
\kappa(t):=\int_0^L k(0,y)\hat{u}(t,y)dy,
\end{equation}
where $\hat u$ is the solution of \eqref{lkdv_ob}.  

Thus we get the following result, which can be compared to \cite{marxcerpa2014outputkdv}. The proof is given in Section \ref{stability_analysis_output_feedback}.\\

\begin{nntheorem}
\label{marxcerpacdc}
For any $\lambda>0$, there exists $C>0$ such that
for any solution of \eqref{lkdv}-\eqref{lkdv_ob}-\eqref{kk-obs} we have
\begin{multline}
\Vert u(t,\cdot)-\hat{u}(t,\cdot)\Vert_{L^2(0,L)}+\Vert \hat{u}(t,\cdot)\Vert_{L^2(0,L)} \leq
C e^{-\lambda t}\Vert u(0,\cdot)\Vert_{L^2(0,L)},\quad \forall t\geq 0.
\end{multline}
\end{nntheorem}

\begin{remark} Notice that from this theorem, we get the exponential decreasing to $0$ of the $L^2$-norm of the solution $u=u(t,x)$. This result is different from \cite{marxcerpa2014outputkdv}, where the initial condition has to be chosen in $H^3(0,L)$. This is due to the fact that the output and the boundary conditions are different.
\end{remark}

\section{Observer design}
\label{observer_design}


The observer \eqref{linear_KdV_observer} is based on a Volterra transformation. It transforms the solution $\tilde{u}:=u-\hat{u}$ which fullfills the following PDE
\begin{equation}
\label{error}
\left\{
\begin{split}
&\tilde{u}_t+\tilde{u}_x+\tilde{u}_{xxx}-p_1(x)[\tilde{u}(t,L)]=0,\\
&\tilde{u}(t,0)=\tilde{u}_x(t,L)=\tilde{u}_{xx}(t,L)=0,\\
&\tilde{u}(0,x)=u_0(x)-\hat u_0(x):=\tilde u_0(x) ,\\
\end{split}
\right.
\end{equation}
into the following PDE
\begin{equation}
\label{error-target}
\left\{
\begin{split}
&\tilde{w}_t+\tilde{w}_x+\tilde{w}_{xxx}+\lambda \tilde{w}=0,\\
&\tilde{w}(t,0)=\tilde{w}_x(t,L)=\tilde{w}_{xx}(t,L)=0,\\
&\tilde{w}(0,x)=\tilde{w}_0(x).
\end{split}
\right.
\end{equation}

We choose the same $\lambda$ than the one used to design the controller. The transformation is given by
\begin{equation}
\label{transformation}
\tilde{u}(x):=\Pi_o(\tilde w(x))=\tilde{w}(x)-\int_x^L p(x,y)\tilde{w}(y)dy, 
\end{equation}
where $p$ is a kernel that satisfies a partial differential equation and will be defined in the following.

In order to find the kernel, we find the relation between \eqref{error} and \eqref{error-target} through the transformation \eqref{transformation}. To do so, we derive the following
\begin{align}
\label{first-kernel}
\tilde{u}_t(t,x)= & \tilde{w}_t(t,y)-\int_x^L p(x,y)\tilde{w}_t(t,y)dy\\
 = & \tilde{w}_t(t,y)+\int_x^L p(x,y)[\tilde{w}_y+\tilde{w}_{yyy}+\lambda \tilde{w}](t,y)dy\nonumber\\
 = & p(x,L)\tilde{w}(t,L)-p(x,x)\tilde{w}(t,x)-p(x,x)\tilde{w}_{xx}(t,x)+p_y(x,x)w_x(t,x)\nonumber\\
&+p_{yy}(x,L)\tilde{w}(t,L)-p_{yy}(x,x)\tilde{w}(t,x)-\int_x^L \tilde{w}(t,y)[p_{yyy}(x,y)+p_y(x,y)-\lambda p(x,y)]dy,
\end{align}where the last line has been obtained performing some integrations by parts.

In addition, we have
\begin{equation}
\label{second-kernel}
\tilde{u}_x(t,x)= \tilde{w}_x(t,x)+p(x,x)\tilde{w}(t,x)-\int_x^L p_x(x,y)\tilde{w}(t,y)dy,
\end{equation}
\begin{equation}
\tilde{u}_{xx}(t,x) =  \tilde{w}_{xx}(t,x)+\frac{d}{dx}p(x,x)\tilde{w}(t,x)+p(x,x)\tilde{w}_x(t,x)
+p_x(x,x)\tilde{w}(t,x)-\int_x^L p_{xx}(x,y)\tilde{w}(t,y)dy,
\end{equation}
\begin{multline}
\label{third-kernel}
\tilde{u}_{xxx}(t,x)= \tilde{w}_{xxx}(t,x)+\frac{d^2}{dx^2}p(x,x)\tilde{w}(t,x)+2\frac{d}{dx}p(x,x)\tilde{w}_x(t,x)
+p(x,x)\tilde{w}_{xx}(t,x)+\frac{d}{dx}p_x(x,x)\tilde{w}(t,x)
\\+p_x(x,x)\tilde{w}_x(t,x)+p_{xx}(x,x)\tilde{w}(t,x)
-\int_x^L p_{xxx}(x,y)\tilde{w}(t,y)dy,
\end{multline}

By adding \eqref{first-kernel}, \eqref{second-kernel} and \eqref{third-kernel}, 
we get
\begin{align}
\tilde{u}_t(t,x)+\tilde{u}_x(t,x)+\tilde{u}_{xxx}(t,x)-p_1(x)\tilde{u}(t,L)=&
\tilde{w}_t(t,x)+\tilde{w}_{x}(t,x)+\tilde{w}_{xxx}(t,x)+\lambda\tilde{w}(t,x)
 \nonumber\\&-\int_x^L \tilde{w}(t,y)[p_{xxx}+p_{yyy}+p_x+p_y-\lambda p](x,y)dy\nonumber\\
&+ \tilde{w}(t,L)[p_{yy}(x,L)+p(x,L)-p_1(x)]\nonumber\\
&+ \tilde{w}(t,x)\left[\frac{d^2}{dx^2}p(x,x)+p_{xx}(x,x)\right.\nonumber\\
&-p_{yy}(x,x)
\left.+\frac{d}{dx}p_x(x,x)-\lambda\right]
\nonumber\\&+\tilde{w}_x(t,x)\left[2\frac{d}{dx}p(x,x)+p_y(x,x)+p_x(x,x)\right].
\end{align}
From this equation, we get the following four conditions.
\begin{itemize}
\item[1.] Equation for all $(x,y)\in\mathcal{T}$:
\begin{equation}
p_{xxx}(x,y)+p_{yyy}(x,y)+p_x(x,y)+p_y(x,y)=
\lambda p(x,y).
\end{equation}
\item[2.] First boundary condition on $(x,x)$ where $x\in [0,L]$:
\begin{equation}
\frac{d^2}{dx^2}p(x,x)+p_{xx}(x,x)-p_{yy}(x,x)+\frac{d}{dx}p_x(x,x)-\lambda=0.
\end{equation}
\item[3.] Second boundary condition on $(x,x)$ where $x\in [0,L]$:
\begin{equation}\label{2nd}
2\frac{d}{dx}p(x,x)+p_y(x,x)+p_x(x,x)=0.
\end{equation}
\item[4.] Appropriate choice of $p_1$:
\begin{equation}
p_1(x)=p_{yy}(x,L)+p(x,L).
\end{equation}
\end{itemize}

Moreover, from the boundary conditions, we get two more restrictions.
\begin{itemize}
\item[5.] To satisfy $\tilde{w}(t,0)=0$, transformation \eqref{transformation} imposes
\begin{equation}
p(0,y)=0,
\end{equation}
which together with \eqref{2nd} implies that 
\begin{equation}
p(x,x)=0.
\end{equation}
\item[6.] To satisfy $\tilde{w}_{xx}(t,L)=0$, transformation \eqref{transformation} imposes
\begin{equation}
p_x(L,L)=0.
\end{equation}
\end{itemize}

Finally, $p$ solves the following equation
\begin{equation}
\label{eq_p}
\left\{
\begin{split}
&p_{xxx}+p_{yyy}+p_{x}+p_y=\lambda p,\quad \forall (x,y)\in\mathcal{T},\\
&p(x,x)=0,\quad \forall x\in [0,L],\\
&p_x(x,x)=\frac{\lambda}{3}(x-L),\quad \forall x\in [0,L],\\
&p(0,y)=0,\quad \forall y\in [0,L].
\end{split}
\right.
\end{equation}
Let us make the following change of variables
\begin{equation}
\bar{x}=L-y,\hspace{0.5cm}\bar{y}=L-x,
\end{equation}
and define $F(\bar{x},\bar{y}):=p(x,y)$. Hence, $F$ satisfies
\begin{equation}
\left\{
\begin{split}
&F_{\bar{x}\bar{x}\bar{x}}(\bar{x},\bar{y})+F_{\bar{y}\bar{y}\bar{y}}(\bar{x},\bar{y})\\
&+F_{\bar{y}}(\bar{x},\bar{y})+F_{\bar{x}}(\bar{x},\bar{y})=-\lambda F(\bar{x},\bar{y})\hspace{0.3cm}(\bar{x},\bar{y})\in\mathcal{T},\\
&F(\bar{x},\bar{x})=0\hspace{0.3cm}\bar{x}\in [0,L],\\
&F_{\bar{x}}(\bar{x},\bar{x})=\frac{\lambda}{3}(L-\bar{x})\hspace{0.3cm}\bar{x}\in [0,L],\\
&F(\bar{x},L)=0\hspace{0.3cm}\bar{y}\in [0,L].
\end{split}
\right.
\end{equation}
This equation has already been studied in \cite{cerpa_coron_backstepping}, where no explicit solution has been found, but where the existence of a solution has been proved. Therefore, we can conclude that  $F=F(x,y)$ and consequently the kernel $p=p(x,y)$ both exist. Note that the function $\Pi_o$ defined by (\ref{transformation}) is linear (by definition) and continuous (because of the existence of $p$). Moreover, $\Pi_o$ is invertible with continuous inverse given by
\begin{equation}
\label{transformation-invP}
\tilde{w}(x)=\Pi^{-1}_o(\tilde{u}(x))=\tilde{u}(x)+\int_x^L m(x,y)\tilde{u}(y)dy
\end{equation}
where $m=m(x,y)$ is also a solution of an equation like \eqref{eq_p} in the triangular domain $\mathcal T$.

\section{Asymptotic stability of the output feedback}
\label{stability_analysis_output_feedback}
The closed-loop system through the transformations \eqref{mapk} and \eqref{transformation} can be written as follows
\begin{equation}
\left\{
\begin{split}
&\hat{w}_t+\hat{w}_x+\hat{w}_{xxx}+\lambda\hat{w}=\\
&-\left\{p_1(x)-\int_x^L k(x,y)p_1(y)dy\right\}\tilde{w}(t,L)\\
&\hat{w}(t,0)=\hat{w}_x(t,L)=\hat{w}_{xx}(t,L)=0\\
&\tilde{w}_t+\tilde{w}_x+\tilde{w}_{xxx}+\lambda \tilde{w}=0\\
&\tilde{w}(t,0)=\tilde{w}_x(t,L)=\tilde{w}_{xx}(t,L)=0.
\end{split}
\right.
\end{equation}

Let us focus on the following Lyapunov function
\begin{equation}
\label{lyapunov_functions}
V(t):=V_1(t)+V_2(t)
\end{equation}
where
\begin{equation}
\label{lyap1}
V_1(t)=A\int_0^L \hat{w}(t,x)^2 dx,
\end{equation}
and
\begin{equation}
\label{lyap2}
V_2(t)=B\int_0^L \tilde{w}(t,x)^2dx.
\end{equation}
The positive values $A$ and $B$ are chosen later. After performing some integrations by parts, we obtain
\begin{equation}
\dot{V}_1(t)\leq \left(-2\lambda+\frac{D^2}{A}\right) V_1(t)+A^2 L |\tilde{w}(t,L)|^2,
\end{equation}
where
\begin{equation}
D=\max_{x\in [0,L]}\left\{p_1(x)-\int_x^L k(x,y)p_1(y)dy\right\}.
\end{equation}
We have also
\begin{equation}
\dot{V}_2(t)\leq -2\lambda V_2(t)-B|\tilde{w}(t,L)|^2.
\end{equation}
Therefore, by choosing
\begin{equation}\label{condA}
A\geq \frac{D^2}{2\lambda}
\end{equation}
and
\begin{equation}\label{condB}
B \geq A^2 L,
\end{equation}
we get
\begin{equation}
\dot{V}(t)\leq -2\mu V(t),
\end{equation}
where
\begin{equation}\label{condmu}
\mu = \left(\lambda-\frac{D^2}{2A}\right)>0.
\end{equation}

This concludes the proof of Theorem \ref{marxcerpacdc} by getting an exponential decay rate equal to $\mu$.

\section{Nonlinear system}
\label{n-l-kdv}
The aim of this section is to prove Theorem \ref{main-th}, i.e., we have to prove the local exponential stability of the nonlinear closed-loop system 
\begin{equation}
\label{nl-outfeed}
\left\{
\begin{split}
&u_t+u_x+u_{xxx}+u_xu=0,\\
&u(t,0)=\kappa(t),\: u_x(t,L)=u_{xx}(t,L)=0,\\
&\hat{u}_t+\hat{u}_x+\hat{u}_{xxx}+\hat{u}_x\hat{u}+p_1(x)[y(t)-\hat{u}(t,L)]=0,\\
&\hat{u}(t,0)=\kappa(t),\: \hat{u}_x(t,L)=\hat{u}_{xx}(t,L)=0,
\end{split}
\right.
\end{equation}
where
\begin{equation}
\kappa(t)=\int_0^L k(0,y)\hat{u}(t,y)dy.
\end{equation}

As before, we consider the evolution of the couple $(\tilde u, \hat u)$ where $\tilde u$ stands for the error $\tilde u=u-\hat u$. Using $\Pi_o$ and its inverse (see \eqref{transformation} and \eqref{transformation-invP}), we define $\tilde w=\Pi_o^{-1}(\tilde u)$. We denote $\hat w=\Pi(\hat u)$, where $\Pi$ is defined in \eqref{mapk}. The inverse $\Pi^{-1}$ is given by 
\begin{equation}
\label{inv-mapk}
\hat{u}(x)=\Pi^{-1}(\hat{w}(x))=\hat{w}(x)+\int_x^L l(x,y)\hat{w}(y)dy,
\end{equation}
where $l$ solves the following equation
\begin{equation}
\left\{
\begin{split}
&l_{xxx}+l_{yyy}+l_x+l_y=\lambda l,\hspace{0.1cm} \text{ in }\mathcal{T},\\
&l(x,L)+l_{yy}(x,L)=0,\hspace{0.1cm} \text{ in }[0,L],\\
&l(x,x)=0,\hspace{0.1cm}\text{ in }[0,L],\\
&l_x(x,x)=\frac{\lambda}{3}(L-x),\hspace{0.1cm}\text{ in }[0,L].
\end{split}
\right.
\end{equation}
The existence of such a kernel $l$ has been proven in \cite[see sections IV and VI]{cerpa_coron_backstepping}. 
Thus, we can see that $(\tilde u, \hat u)$ is mapped into $(\tilde w, \hat w)$ solution of the coupled target system
\begin{equation}
\begin{split}
\hat{w}_t(t,&x)+\hat{w}_x(t,x)+\hat{w}_{xxx}(t,x)+\lambda\hat{w}(t,x)\\
=&-\left\{p_1(x)-\int_x^L k(x,y)p_1(y)dy\right\}\tilde{w}(t,L)\\
&-\left(\hat{w}(t,x)+\int_x^L l(x,y)\hat{w}(t,y)dy\right)\\
&\cdot\left(\hat{w}_x(t,x)+\int_x^L l_x(x,y)\hat{w}(t,y)dy\right)\\
&-\frac 1 2 \int_x^L| \hat u(t,y)|^2 k_y(x,y)dy,
\end{split}
\end{equation} 
\begin{equation}\label{eq_wtildeNL}
\begin{split}
\tilde{w}_t(t,&x)+\tilde{w}_x(t,x)+\tilde{w}_{xxx}(t,x)+\lambda\tilde{w}(t,x)\\
=&-\left(\tilde{w}(t,x)-\int_x^L p(x,y)\tilde{w}(t,y)dy\right)\\
&\cdot\left(\tilde{w}_x(t,x)-\int_x^L p_x(x,y)\tilde{w}(t,y)dy\right)\\
&-\left(\hat{w}(t,x)+\int_x^L l(x,y)\hat{w}(t,y)dy\right)\\
&\cdot\left(\tilde{w}_x(t,x)-\int_x^L p_x(x,y)\tilde{w}(t,y)dy\right)\\
&-\left(\tilde{w}(t,x)-\int_x^L p(x,y)\tilde{w}(t,y)dy\right)\\
&\cdot\left(\hat{w}_x(t,x)+\int_x^L l_x(x,y)\hat{w}(t,y)dy\right)\\
&+ \int_x^L\Big[\frac{| \tilde u(t,y)|^2}{2}+\tilde u(t,y)\hat u(t,y) \Big]m_y(x,y)dy,
\end{split}
\end{equation}
with boundary conditions
\begin{equation}
\hat{w}(t,0)=\hat{w}_x(t,L)=\hat{w}_{xx}(t,L)=0,
\end{equation}
\begin{equation}
\tilde{w}(t,0)=\tilde{w}_x(t,L)=\tilde{w}_{xx}(t,L)=0.
\end{equation}

As in previous section, we will prove the stability of this system by using the same Lyapunov function (\ref{lyapunov_functions}). We derivate (\ref{lyap1}) with respect to time as follows
\begin{equation}\label{pre1}
\dot{V}_1(t) =2A\int_0^L \hat{w}_t(t,x)\hat{w}(t,x)dx
\leq \left(-2\lambda+\frac{D^2}{A}\right)V_1(t)+A^2 L |\tilde w(t,L)|^2
-2A\int_0^L\hat{w}(t,x)F(t,x)dx,
\end{equation}
where
\begin{multline}
F(t,x)= \hat{w}(t,x)\hat{w}_x(t,x)+\hat{w}(t,x)\int_x^L l_x(x,y)\hat{w}(t,y)dy
+\hat{w}_x(t,x)\int_x^L l(x,y)\hat{w}(t,y)dy\\
+\left(\int_x^L l(x,y)\hat{w}(t,y)dy\right)\left(\int_x^L l_x(x,y)\hat{w}(t,y)dy\right)
+\frac 1 2 \int_x^L| \hat u(t,y)|^2 k_y(x,y)dy.
\end{multline}

By using the same argument as in \cite{cerpa_coron_backstepping,coron-lu}, we can prove the existence of a positive constant $K_1=K_1(\Vert l\Vert_{C^1(\mathcal T)},\Vert k\Vert_{C^1(\mathcal T)})$ such that
\begin{equation}
\left|A\int_0^L \hat{w}(t,x)F(t,x)dx\right|\leq K_1\left(\int_0^L |\hat{w}(t,x)|^2dx\right)^{\frac{3}{2}}.
\end{equation}

Then, we estimate $\dot V_2(t)$ as follows
\begin{equation}\label{pre2}
\dot{V}_2(t)\leq -2\lambda V_2(t)-2B\int_0^L \tilde{w}(t,x)G(t,x)dx-B|\tilde{w}(t,L)|^2
\end{equation}
where $G=G(t,x)$ is the right-hand side of \eqref{eq_wtildeNL}. As before, we can prove the existence of a positive constant $K_2=K_2(\Vert l\Vert_{C^1(\tau)},\Vert p\Vert_{C^1(\tau)},\Vert m\Vert_{C^1(\tau)})$ such that
\begin{equation}
\left|2B\int_0^L \tilde{w}(t,x)G(t,x)dx\right|\leq K_2\left(\int_0^L |\hat{w}(t,x)|^2dx\right)^{\frac{3}{2}}
+K_2\left(\int_0^L |\tilde{w}(t,x)|^2dx\right)^{\frac{3}{2}}.
\end{equation}

Therefore, for $A,B,\mu$ satisfying \eqref{condA}-\eqref{condB}-\eqref{condmu}, we have
\begin{equation}
\begin{split}
\dot{V}(t)\leq &-2\mu V(t)+K_1\left(\int_0^L |\hat{w}(t,x)|^2dx\right)^{\frac{3}{2}}\\
&+K_2\left(\int_0^L |\hat{w}(t,x)|^2dx\right)^{\frac{3}{2}}
\\
&+K_2\left(\int_0^L |\tilde{w}(t,x)|^2dx\right)^{\frac{3}{2}}.
\end{split}
\end{equation}
If there exists $t_0\geq 0$ such that
\begin{equation}
\Vert \tilde{w}(t_0,.)\Vert_{L^2(0,L)}\leq \frac{\mu}{K_2}
\end{equation}
and
\begin{equation}
\Vert \hat{w}(t_0,.)\Vert_{L^2(0,L)}\leq \frac{\mu}{K_1+K_2}
\end{equation}
we can conclude
\begin{equation}
\dot{V}(t)\leq -\mu V(t),\quad \forall t\geq t_0.
\end{equation}
Thus, we get
\begin{equation}
\Vert \tilde{w}(t,.)\Vert_{L^2(0,L)}+\Vert \hat{w}(t,.)\Vert_{L^2(0,L)}\leq  e^{-\frac{\mu}{2}t}\left(\Vert \tilde{w}_0\Vert_{L^2(0,L)}+\Vert \hat{w}_0\Vert_{L^2(0,L)}\right),\quad \forall t\geq 0,
\end{equation}
provided that
\begin{equation}
\begin{split}
\Vert \hat{w}_0\Vert_{L^2(0,L)}\leq & \frac{\mu}{K_1+K_2},\\
\Vert \tilde{w}_0\Vert_{L^2(0,L)}\leq & \frac{\mu}{K_2}.
\end{split}
\end{equation}
This concludes the proof of Theorem \ref{main-th} by getting the exponential decay of the system with a smallness condition on the $L^2$-norm of the initial data $u_0,\hat u_0$.

\section{Numerical Simulations}
\label{simu}

In this section we provide some numerical simulations showing the effectiveness of our control design. In order to discretize our KdV equation, we use a finite difference scheme inspired from \cite{nm_KdV}. The final time for simulations is denoted by $T_{final}$. We choose $(N_x+1)$ points to build a uniform spatial discretization of the interval $[0,L]$ and $(N_t+1)$ points to build a uniform time discretization of the interval $[0,T_{final}]$. Thus, the space step is $\Delta x=L/N_x$ and the time step $\Delta t=T_{final}/N_t$. We approximate the solution with the notation $u(t,x)\approx U^i_j$, where $i$ and $j$ refer to time and space discrete variables, respectively.

Some used approximations of the derivative are given by
\begin{equation}
u_x(t,x)\approx \nabla_-(U^i_j)=\frac{U^i_j-U^i_{j-1}}{\Delta x}
\end{equation}
or
\begin{equation}
u_x(t,x)\approx \nabla_+(U^i_j)=\frac{U^i_{j+1}-U^i_{j}}{\Delta x}.
\end{equation}
As in \cite{nm_KdV}, we choose rather the following
\begin{equation}
u_x(t,x)\approx\frac{1}{2}(\nabla_++\nabla_-)(U^i_j)=\nabla(U^i_j).
\end{equation}
For the other differentiation operator, we use
\begin{equation}
u_{xxx}(t,x)\approx \nabla_+\nabla_+\nabla_-(U^i_j)
\end{equation}
and
\begin{equation}
u_t(t,x)\approx \frac{U^{i+1}_j-U^i_j}{\Delta t}.
\end{equation}
%
Let us introduce a matrix notation. Let us
consider $D_-,D_+,D \in\mathbb{R}^{N_x\times N_x}$  given by

\begin{equation}
D_-=\frac{1}{\Delta x}\begin{bmatrix}
1 & 0 & \ldots & \ldots & \,\,0\\
-1 & 1 & \ddots & &  \,\,\vdots\\
0 & -1 & \ddots & \ddots &  \,\,\vdots\\
\vdots & \ddots & \ddots & 1 &  \,\,0\\
0 & \ldots & 0 & -1 &  \,\,1
\end{bmatrix}, 
\end{equation} 
\begin{equation}
D_+=\frac{1}{\Delta x}\begin{bmatrix}
-1 & 1 & 0 & \ldots & 0\\
0 & -1 & 1 & \ddots & \vdots\\
\vdots & \ddots & \ddots & \ddots & 0 \\
\vdots &  & 0 & -1 & 1\\
0 & \ldots &\ldots & 0 & -1
\end{bmatrix},
\end{equation}
and $D=(D_++D_-)/2$. Let us define $A:=D_+D_+D_-+D$, and $C:=A+\Delta t I_{d}$ where $I_{d}$ is the identity matrix.  Moreover, we will denote, for each discrete time $i$, $$U^i:=\begin{bmatrix}
U^i_1 & U^i_2 & \ldots & U_{N_x+1}^i
\end{bmatrix}^T$$ the plant state, and $$O^i:=\begin{bmatrix}
O^i_1 & O^i_2 & \ldots & O_{N_x+1}^i
\end{bmatrix}^T$$ the observer state. The state output will be denoted by $Y_U$ and 
$Y_O$ stands for the observer output.  The discretized controller gain $K$ and observer gain $P$, respectively, are defined by
$$K=\begin{bmatrix}
K_1 & K_2 & \ldots & K_{N_x+1}
\end{bmatrix}^T$$ and $$P=\begin{bmatrix}
P_1 & P_2 & \ldots & P_{N_x+1}
\end{bmatrix}^T.$$ We compute them from a successive approximations method (see \cite{cerpa_coron_backstepping}). 
Since we have the nonlinearity $uu_x$, we use an iterative fixed point method to solve the nonlinear system $$CU^{i+1}=U^i-\frac{1}{2}D(U^{i+1})^2.$$ With $N_{iter}=5$, which denotes the number of iterations of the fixed point method, we get good approximations of the solutions.

Given $U^0$, $O^0=0$, $K$, and $P$, the following is the structure of the algorithm used in our simulations. \\

\fbox{\begin{minipage}{16cm}
\textbf{While $i < N_t$}\\
\begin{itemize}
\item[$\bullet$] $U_{1}^{i+1}=O_1^{i+1}=\sum_{j=1}^{N_x} \Delta x\frac{O^i_{j+1}K_{j+1}+O^i_jK_j}{2}$\\
$U^{i+1}_{N_x}=U^{i+1}_{N_{x}+1}=U^{i+1}_{N_x-1},\: O^{i+1}_{N_x}=O^{i+1}_{N_x+1}=O^{i+1}_{N_x+1}$;\\
\item[$\bullet$] $Y_U:=U^i(N_x+1)$, $Y_O:=O^i(N_x+1)$;\\
\item[$\bullet$] Setting $J(1)=O^i$, for all $k\in\lbrace 1,\ldots, N_{iter}\rbrace$, solve \\
$$J(k+1)=C^{-1}(O^i-\frac{1}{2}D(J(k))^2+P(Y_U-Y_O))$$\\
Set $O^{i+1}=J(N_{iter})$;\\
\item[$\bullet$] Setting $\tilde{J}(1)=U^i$, for all $k\in\lbrace 1,\ldots, N_{iter}\rbrace$, solve \\
$$\tilde{J}(k+1)=C^{-1}(U^i-\frac{1}{2}D(\tilde{J}(k))^2)$$\\
Set $U^{i+1}=J(N_{iter})$;\\
\item[$\bullet$] $t=t+dt$;\\
\end{itemize}
\textbf{End}
\end{minipage}}
\vspace{0.3cm}

\begin{figure}[h!]
\centering{
\includegraphics[scale=0.6]{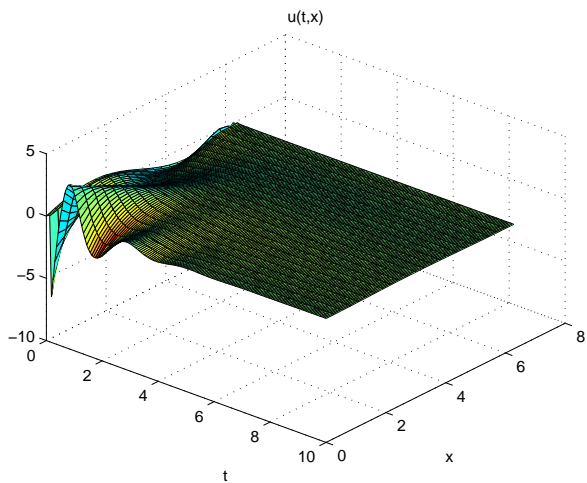}
\caption{Solution of the closed-loop system \eqref{lKdV_equation_ld}-\eqref{measurement}-\eqref{linear_KdV_observer}-\eqref{kk-obs}}
\label{sol}
}
\end{figure}

\begin{figure}[h!]
\centering{
\includegraphics[scale=0.6]{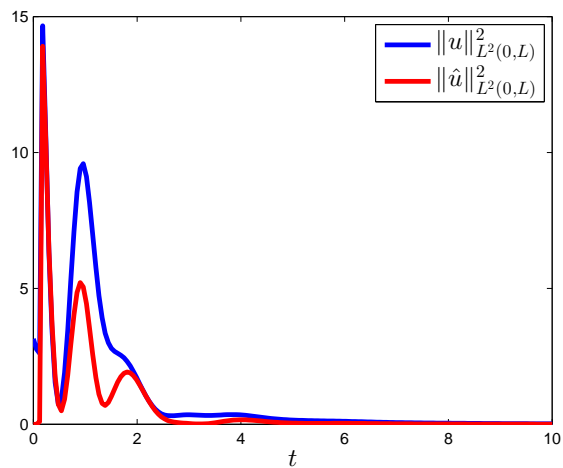}
\caption{Time evolution of the $L^2$-norm for the state (blue line) and the observer (red line).}
\label{lyap-figure}
}
\end{figure}

\begin{figure}[h!]
\centering{
\includegraphics[scale=0.6]{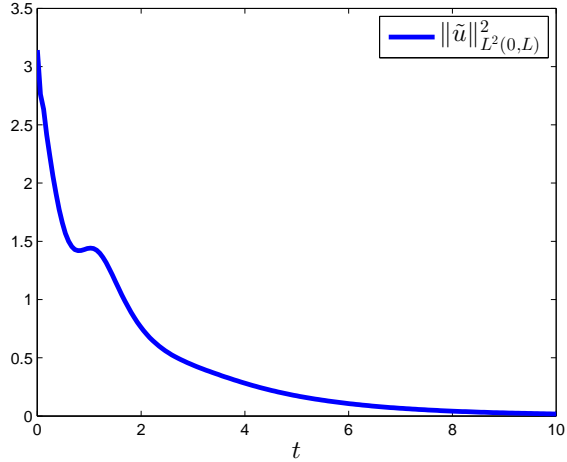}
\caption{Time evolution of the $L^2$-norm for the observation error $u-\hat{u}$.}
\label{lyap-error-figure}
}
\end{figure}

In order to illustrate our theoretical results, we perform some simulations on the domain $[0, 2\pi]$. We take $N_x=30$, $N_t=167$, $T_{final}=10$, $\lambda=2$, $u_0(x)=sin(x)$ and $\hat u_0(x)=0$. Figure \ref{sol} illustrates the convergence  to the origin of the solution of the closed-loop system \eqref{lKdV_equation_ld}-\eqref{measurement}-\eqref{linear_KdV_observer}-\eqref{kk-obs}. Figure \ref{lyap-figure} illustrates the $L^2$-norm of this solution and the $L^2$-norm of the solution of the observer \eqref{linear_KdV_observer}. Finally, Figure \ref{lyap-error-figure} shows that the observation error $(u-\hat{u})$ converges to $0$  in $L^2$-norm. From the simulations, this convergence seems to be exponential as expected.

\section{Conclusion}
\label{con}

In this paper, an output feedback control has been designed for the Korteweg-de Vries equation. This controller uses an observer and gives the local exponential stability of the closed-loop system. Numerical simulations have been provided to illustrate the efficiency of the output feedback design.

In order to go to a global feedback control for the nonlinear KdV equation, a first step should be to build some nonlinear boundary controls giving a semi-global exponential stability. That means that for any fixed $r>0$ we can find a feedback law exponentially stabilizing to the origin any solution with initial data $u_0$ whose $L^2$-norm is smaller than $r$. The term semi-global comes from the fact that the decay rate can depend on $r$. Some internal feedback controls are given for KdV in \cite{pazoto}, \cite{rosier-zhang2} and \cite{mcpa-siam}. The latter considers saturated controls.

The second step, in order to deal with the output case, is to design an observer for the nonlinear KdV equation to obtain a global asymptotic stability. We could for instance follow the same approach than in the finite-dimensional case performed in \cite{GHO}.


\bibliographystyle{plain}
\bibliography{Biblio}

\end{document}